\documentclass[12pt,reqno]{amsart}
\usepackage{hyperref}
\usepackage{amsmath,amssymb,amsthm,amsfonts}
 \usepackage{bm,latexsym,commath} 

\usepackage{graphicx,float,multicol}

\usepackage{siunitx,caption}
\usepackage{array, booktabs} 

\begin{document}

\title{Exact solutions to the cancer laser ablation modeling}

\author{Luisa Consiglieri}
\address{Luisa Consiglieri, Independent Researcher Professor, European Union}
\urladdr{\href{https://sites.google.com/site/luisaconsiglieri}{https://sites.google.com/site/luisaconsiglieri}}

\begin{abstract} 
The present paper deals with the study of the fluence rate over both healthy and tumor tissues in the presence of focal laser ablation (FLA).
We propose new analytical solutions for the coupled partial differential equations (PDE) system,
which includes the transport equation modeling the light penetration into biological tissue, the bioheat equation modeling the heat transfer and its respective damage.
The present building could be the first step to the knowledge of the mathematical fr
amework for biothermophysical problems, as well as the main key to simplify the numerical calculation due to its no cost.
We derive exact solutions and simulate results from them.
We discuss the potential physical contributions and present respective conclusions about
 (1) the validness of  the diffusion  approximation of the radiative transfer equation; 
(2) the local behavior of the source of scattered photons; 
(3) the unsteady-state of the fluence rate; and (4) the boundedness of the critical time of the thermal damage to the cancerous tissue.
We also discuss some controversial and diverging hypotheses.
\end{abstract}


\maketitle

\section{Introduction}
\label{intro}

Thermal ablation is a widely accepted therapy for denaturing proteins and then the destruction of tissue cells.
It has been used both for the treatment of atrial fibrillation \cite{voel} or  cancer \cite{fan,Loiola}.
In the last decades, several energy sources have been applied (see \cite{Laza,HeBischof} and references therein), including radiofrequency \cite{consiglieri2003,lc2012,pupo,g-suarez},
 microwave \cite{chiang,WangWuWu}, cryoablation \cite{Franco,merryman,kumar}, ultrasound  \cite{ho}, and laser \cite{mordon,vogel}.
Focal Laser Ablation (FLA) is a minimally invasive procedure used to treat tumors by using laser energy to heat and destroy cancerous cells.
Indeed, FLA is gaining acceptance as a safe option and effective alternative to surgical resection \cite{putzer,DeVita,schena}
and the references therein.

An ideal hyperthermia treatment is  one that selectively destroys the diseased cells without damaging the surrounding healthy tissue. 
For this reason, it is still essential  build analytical solutions to respond such selective targeting.
The FLA model results from three distinct phenomena: conversion of laser light into heat, transfer of
heat, and coagulation necrosis. It consists of  the coupled partial differential equations (PDE) system which includes
 the  diffusion approximation of the radiative transfer equation modeling the light penetration into biological
tissue,  the Pennes bioheat transfer equation modeling the heat transfer, and the Arrhenius burn integration modeling the dimensionless indicator of damage.
The worldwide literature of the heat transfer on the ablation modeling commonly considers the heat source as constant or of particular profile.
In FLA, the heat source, the so-called absorbed optical power density,  is induced 
by the conversion of laser light into heat and has a well-known linear relationship with the fluence rate.  
In the presence of this relationship,  accurate knowledge of the fluence rate  is the  key feature for the mathematical framework.

For a light source acting as a cylindrical diffuser, the fluence rate is obtained
in terms of zero-order modified Bessel functions \cite[pp. 172-174]{w-vang6}.
Also other expressions exist (see, for instance, \cite{liu-boas} and the references therein).
For a single spherically symmetric point source emitting,  we mention
\cite{oden,szhu,w-vang6} and references therein.
The laser beam has  a crater-like shape, since the intensity is not constant across the beam width but decreases
from the center towards the edge \cite{whit1d}.

Some analytical solutions to the bioheat transfer problem are studied for RF ablation in  multiregion domains
\cite{durk,lc2013,lc2015}. We refer to \cite{oden} in which  finite element model
simulations of laser therapy are computed.
For  the skin surface in \cite{liu}, the authors  provide a theoretical study of the thermal wave effects
in the bioheat transfer problem involving high heat flux incident with a short duration.
For breast tumor thermal therapy in \cite{ho},
the authors propose cylindrical ultrasound phased array with a multifocus pattern scanning strategy.
 For  a spherical tumor  in \cite{pupo}, the authors simulate the current density distribution 
when a 3D current density is generated by different multiple-straight needle electrode configurations.
By using finite element modeling in \cite{paul}, the authors analyze
the photo-thermal heating with near-infrared radiation in the presence of intravenous blood injection or
intratumorally injected gold nanorods.
We refer to  \cite{lc2024} the study of the endovenous laser ablation (EVLA), which stands for 
 a similar technique  to treat varicose veins by using laser energy to close off the affected veins.

One-dimensional, time-dependent models for the removal of
material using a laser beam were published in \cite{whit1d} and the references therein. 
The thermal penetration differs with the laser system (Argon, Nd:YAG, and CO\(_2\) \cite{mordon}). 
Failure rates may depend on the  design of the optical fibers. 
In \cite{bane}, the authors  show that
the heat-affected zone is significantly reduced by using a short pulsed
laser of the same average power (\SI{150}{\milli\watt}  for 25 seconds) as compared to a CW laser source. 
Short pulsed lasers are being used in a variety of applications
such as remote sensing, optical tomography, laser surgery, and ablation processes.

In the present work, we derive exact solutions and apply the resulting solutions under Nd:YAG and diode  lasers for their wide application.
We analyze the light source over different tissues and its tissue absorption under the Beer--Lambert  (BL) law.
We focus our attention on breast and prostate tumors that are common cancers worldwide diagnosed in women and men, respectively \cite{Sfarra,manenti}.

The outline of the present work is as follows.
We firstly describe the mathematical model in Section \ref{model}. In Section \ref{sanalytical}, we
 propose new exact solutions for the unknown functions under study. In Section \ref{results}, 
we apply the previous results to the fluence rate and temperature in the existence sections \ref{fluencerate} and \ref{temperature},  respectively,
and present some simulation  in Section \ref{profiles}. 
Section \ref{discuss} indicates the main interpretations and conclusions.
Finally, two Appendices are included to contain mathematical proofs relative to  Section \ref{results}.

\section{Mathematical formulation}
\label{model}

\begin{figure}[H]
\includegraphics[width=14 cm]{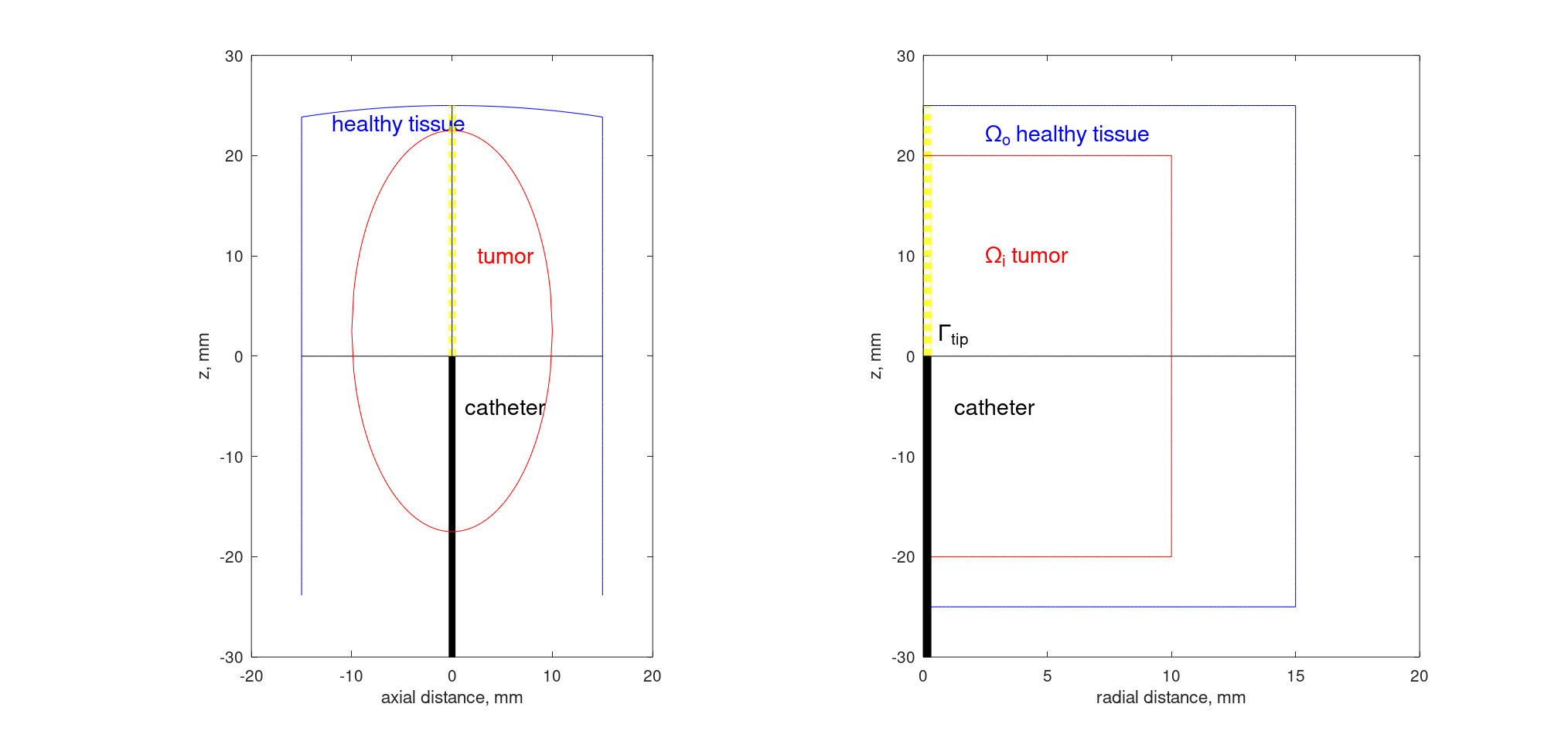}
\caption{\footnotesize Schematic representations of the laser-tissue system: the laser emission tip
\( \Gamma_\mathrm{tip} =\{(x,y)\in \mathbb R^2 :\,\,x^2+y^2<r_{\mathrm{f}}^2 \} \times \{0\} \),
the light-emitting cylinder \(]0;r_\mathrm{f}[\times ]0;L[\),
 the tumor  \( \Omega_\mathrm{i} =]0; r_\mathrm{i}[\times ] -\ell ;\ell[ \),
and the surrounding healthy tissue  \( \Omega_\mathrm{o}= \{(x,y)\in\mathbb R^2:\,\ x^2+y^2<r_\mathrm{i} ^2 \}\times ]\ell;L[\cup
 \{(x,y)\in\mathbb R^2:\,\, r_\mathrm{i}^2<x^2+y^2<r_\mathrm{o} ^2 \}\times ]0;L[\).
Left: Sagittal view. Right: Cylindrical coordinates \((r,z)\).\label{domain}}
\end{figure}

We assume the geometry of the laser--tissue system to be as illustrated in Figure~\ref{domain}.
The laser applicator usually consists of an optical fiber  (with radius \(r_\mathrm{f}\)) 
incorporated within  a catheter.
The diameter of the optical fiber for medical applications varies between  \SI{200}{\micro\metre} and
 \SI{600}{\micro\metre} \cite{niemz},  while the laser beam diameter is \SI{2}{\milli\metre} \cite{bane}.
 Thus, we may assume  \(r_\mathrm{f}=\SI{0.25}{\milli\metre}\)  \cite{marq}.

Both breast and prostate tumors can be modeled as ellipsoid (see \cite{Sfarra,marq} and references therein).
We assume
the tumor \(\Omega_\mathrm{i}\) and surrounding healthy  tissue \(\Omega_\mathrm{o}\)  to be cylindrical with
\(r_\mathrm{i}\) and \(r_\mathrm{o}\) being the inner and outer radius, respectively.
Let us define the bidomain
\[\Omega =\mbox{int} (\overline \Omega_\mathrm{i}\cup \Omega_\mathrm{o}).\]
Mean breast and prostatic carcinoma sizes may vary in different ranges.  We assume the inner radius  \(r_\mathrm{i}=\SI{1}{\centi\metre}\).
Although the primary goal is to destroy all malignant tissue, the surrounding nonmalignant tissue should be preserved. To evaluate this rim  of surrounding tissue,
the outer radius  \(r_\mathrm{o}\) will be determined such that the fluence rate vanishes.
Indeed, the respective photons just "vanish" from the system due to the light reaching the system border is absorbed.

The laser light is emitted from the distal end
of the fiber on an active length which is approximately half the tumor size \cite[p. 82]{niemz}.
By this reason,  the distal end is assumed to be centered in the middle of the tumor, that is, the longitudinal coordinate is \(-\ell<z<\ell\).       

Three different instants of time are of reference:
 the temporal pulse width \(t_p\)  (also known as pulse duration or pulse  width, for short), the  critical time \(t_\mathrm{crit}\),  and the exposure duration \(t_\mathrm{end}\).
The pulse width  \(t_p\) is the time measured across a pulse,  at its full width half maximum (FWHM).
The  critical time \(t_\mathrm{crit}\) is introduced in Section \ref{sdam}.
The exposure duration \(t_\mathrm{end}\) consists of \(N\) pulse widths 
 and \(N-1\) pulse-to-pulse intervals  \(\Delta t\) (the so-called periods), which are ten times the pulse width \cite{mordon}.

\subsection{Photon transport}
\label{s-radiative}

The photon transport is governed by the diffusion approximation of the radiative transfer equation (see, for instance, \cite{niemz,w-vang6,w-vang7})
\begin{equation}\label{light}
\frac{1}{\nu}\frac{\partial\phi}{\partial t}-D\Delta \phi +\mu_\mathrm{a}\phi=S
\qquad \mathrm{in }\ \Omega\times ]0 ; t_\mathrm{end} [,
\end{equation}
where:
\begin{itemize}
\item \(\phi\) is the fluence rate [\si{\watt\per\square\metre}];
\item  \(\nu\)  is the  light velocity in the tissue  [\si{\metre\per\second}], \(\nu=c/n\), with
 the  light velocity   \(c=\SI{3e8}{\metre\per\second}\) and the relative refractive index  \(n\);
\item  \(\mu_\mathrm{a}\) is the absorption coefficient [\si{\per\metre}].
\end{itemize}

The diffusion coefficient is 
\[ D=\frac{1}{3(\mu_\mathrm{a}+\mu'_\mathrm{s})}. \] 
 Here,  the reduced scattering coefficient \(\mu'_\mathrm{s}\)  [\si{\per\metre}] is a function of the wavelength  \(\lambda\),
which is normalized by the reference wavelength \SI{500}{\nano\metre}, 
\[
\mu'_\mathrm{s}=a\left( \frac{\lambda}{500} \right) ^{-b}
\]
where \(a\) and \(b\) are known constants \cite{jacq}. 
The reduced scattering coefficient  obeys \(\mu'_\mathrm{s}=(1-g)\mu_\mathrm{s}\), with
\(g\) being the scattering anisotropy coefficient and \(\mu_\mathrm{s}\) being the scattering coefficient.
The biological tissue is a strongly scattering media and then \(g\) ranges from \(0.7\) to \(0.99\) for most biological tissues \cite{niemz}.
The optical parameters of the tissue are known (cf. Table~\ref{tabopt}), which are assumed to be  temperature independent  \cite{whit1d}.
The average of \(\mu_\mathrm{a}\) and \(\mu'_\mathrm{s}\) belong to the range 0--0.6 and 0--20 \si{\per\centi\metre},  respectively.
\begin{table}[H]
\caption{Average optical parameters for breast and prostate (tumor and healthy) tissues \cite{jacq,niemz,szhu}. \label{tabopt}}       
\begin{tabular}{llllllll}
\toprule
 && \(\mu_\textrm{a}\)  &  [\si{\per\centi\metre} ] & & \(\mathbf{a}\)  [\si{\per\milli\metre} ] & \(\mathbf{b}\) & \(\mathbf{n}\) \\
\midrule
\(\lambda\)  [\si{\nano\metre}] & 810  & 980 &1,064& && &\\ 
\midrule
Breast tumor &0.08& 0.07 & 0.06& & 2.07 & 1.487 &1.4 \\
Prostate tumor & 0.12 & 0.11 & 0.1 &  & 3.36 & 1.712 &1.4 \\
Breast tissue &0.17 &0.2 & 0.3 & & 1.68& 1.055 &1.35 \\
Prostate tissue &0.6 & 0.5& 0.4& & 3.01 & 1.549 & 1.37 \\
\bottomrule
\end{tabular}
\end{table}

We consider that the beam divergence is negligible around the focus, and then
 the effects of the laser radiation in terms of the
local absorption of light are according to the Beer--Lambert law   \cite{w-vang7}.
The source of scattered photons  \(S\)  [\si{\watt\per\cubic\metre}] represents the power injected in the unit volume, and it is given by
\begin{equation}
S(r,z,t)= \frac{\mu_\mathrm{s}(\mu_\mathrm{t}+g\mu_\mathrm{a})}{\mu_\mathrm{a}+ \mu_\mathrm{s} '} E(r,t)  \exp[-\mu_\mathrm{t}z] ,
\label{source}
\end{equation}
for \((r,z,t)\in ]0; r_\mathrm{o}[ \times ]0; L[\times ] 0; t_\mathrm{end} [\).
Here,  the total attenuation coefficient is
\begin{equation}\label{mut}
\mu_\mathrm{t}=\mu_\mathrm{a}+\mu_\mathrm{s} = \mu_\mathrm{a}+(1-g)^{-1} a(500/\lambda)^b,
\end{equation}
and  the planar irradiance,  over the exposure duration \(t_\mathrm{end}\),  is
\begin{equation}\label{planar}
E (r,t)=  \frac{P_\mathrm{peak}  }{\pi r_\mathrm{f}^2} \chi_{[0;r_\mathrm{f} ]} (r) \sum_{j=0}^{N-1}\chi_{[t_j;t_j+ t_p ]} (t) ,
\end{equation}
with \( P_\mathrm{peak} \)  [\si{\watt}] standing for the maximum optical power output by the laser,  each pulse time \(t_j = j( t_p +\Delta t)\), 
 and \(\chi_{I} \) standing   for the characteristic function over the interval \(I\).
The number of pulses is such that \(N= ( t_\mathrm{end}+\Delta t)/( t_p +\Delta t)\).

The difference on relative refractive indices  between tumor media and healthy tissue is positive,  namely
the tumor refractive index \(n_0 =1.4>n_1\) (the refractive index of the medium of the transmitted  ray).
At the interfaces,  where mismatched refractive indices occur,  we may consider the Robin boundary conditions \cite{Capart,w-vang7}
\begin{align}
-D\frac{\partial \phi}{\partial z} +\gamma_r \phi = - S_1 \quad &\mbox{ at } z =\ell; \label{BCz}\\
-2D\frac{\partial \phi}{\partial r} +\gamma_r \phi = 0 \quad &\mbox{ at } r = r_\mathrm{i}.\label{BCr}
\end{align}
Here, the coefficient \(\gamma_r\) denotes a reflectance dependent function \cite{w-vang6,w-vang7}, 
and \(S_1(r,z)= g\mu_\mathrm{s}/ \mu_\mathrm{tr}*  E(r,t_p)  \exp[-\mu_\mathrm{t}z] \),  with \(\mu_\mathrm{tr}=\mu_\mathrm{a}+\mu_\mathrm{s}'\) being 
  the transport attenuation coefficient.
Notice that the source of scattered photons \(S = -\mbox{div} S_1 +S_2\), 
where \(S_2(r,z)= \mu_\mathrm{s}*  E(r,t_p)  \exp[-\mu_\mathrm{t}z] \)  and 
\begin{equation}\label{defS1}
S_1(r,z)= \frac{g}{\mu_\mathrm{t}+g\mu_\mathrm{a}} S(r,z,t_p) .\end{equation}
The Robin boundary conditions \eqref{BCz}-\eqref{BCr} mean whenever  light reaching the system borders 
some photons are reflected back into the system.

The remaining interfaces, \(z=0\) and 
\(r=r_\mathrm{f}\), are indeed  mathematical boundaries and not tissue borders. Thus, we assume continuity of the fluence rate \(\phi\) and its flux.

\subsection{Heat transfer}
\label{s-heat}

The heat transfer due to the energy of light deposited is governed by the  Pennes bioheat transfer equation that distinguishes  from nonliving medium:
\begin{equation}
\rho c_\mathrm{p}\frac{\partial T}{\partial t} - \nabla \cdot (k\nabla T) 
 + c_\mathrm{b} \omega_\mathrm{b}  (T-T_\mathrm{b}) =  q  \qquad \mathrm{in } \ \Omega\times ]0 ; t_\mathrm{end}[,
\label{bioheat}
\end{equation}
where:
\begin{itemize}
\item \(T\) is the temperature [\si{\kelvin}];
\item \(\rho\) is the density of the tissue  [\si{\kilogram\per\cubic\metre}];
\item \(\omega_\mathrm{b} \) represents the blood perfusion rate  [\si{\kilogram\per\cubic\metre\per\second}]   that occurs in the capillary bed.
\end{itemize}
In particular,  the constant \(T_\mathrm{b}\) represents the temperature  of the blood (assumed to be \SI{38}{\degreeCelsius}),
\(\rho_\mathrm{b}\)  and \(c_\mathrm{b}\) denote the density  and the specific heat capacity of the blood, respectively, and
 \(c_\mathrm{b} \omega_\mathrm{b} \) accounts  for the heat conducted in direction of the
contribution of flowing blood to the overall energy balance, before the critical coagulation temperature.
The blood perfusion rate \(\omega_\mathrm{b} \) may have a nonlinear temperature dependence \cite{ryl2005}. 
Here, we assume  \(\omega_\mathrm{b} \) that obeys the linear time dependence \eqref{omega}.

The specific heat capacity per unit mass  \( c_\mathrm{p}\)  [\si{\joule\per\kilogram\per\kelvin}]
and thermal conductivity  \(k\)  [\si{\watt\per\metre\per\kelvin}] may be based on the following relationships \cite[pp. 68-69]{niemz}:
\begin{align*}
c_\mathrm{p}&=1550+2800\frac{\rho_\mathrm{w}}{\rho} ; \\ 
k&=0.06+0.57 \frac{\rho_\mathrm{w}}{\rho} , 
\end{align*}
where \(\rho_\mathrm{w}\) denotes the density of the water.

The metabolism behavior is neglected due to its small contribution to the temperature response.  
Then, the laser light induced heat source \(q\)  [\si{\watt\per\cubic\metre}],
the so-called absorbed optical power density,  is governed by
\begin{equation}\label{power}
q=\mu_\mathrm{a} \phi.
\end{equation}

We refer to \cite{guo-kim,jaun,tang} the thermal relaxation time in the conduction equation with both Fourier and non-Fourier effects,
instead of considering the system \eqref{bioheat}-\eqref{power}.
In the last years, the researchers keep proposing various alternatives to the benchmark work developed by Pennes.
We may mention \cite{dutta}, and the references therein,  in which a local thermal non-equilibrium bioheat model is studied.

\subsection{Thermal damage to the tissue}
\label{sdam}

Thermal laser-tissue interaction
exhibits four phases: coagulation, vaporization, carbonization, and melting (vacuolation) \cite[pp. 58-63]{niemz}.

A critical time \(t_\mathrm{crit}\), which corresponds to
  the dimensionless indicator of damage when it is equal to one: \(\Omega(t_\mathrm{crit} )=1\), 
 obeys  the Arrhenius burn integration \cite{ryl2005}
\begin{equation}\label{arr}
\frac{1}{A}=\int_{0}^{t_\mathrm{crit}} \exp\left[-\frac{E_\mathrm{a} }{ RT(r,z,\tau)}\right]\dif{\tau},
\end{equation}
where \(R\) is the universal gas constant (\SI{8.314}{\joule\per\mol\per\kelvin}),
\( A\) is the frequency factor [\si{\per\second}] (also known as Arrhenius factor), and 
\(E_\mathrm{a}\) is the activation energy for the irreversible damage reaction [\si{\joule\per\mol}].

The linear behavior of the  blood perfusion rate \(\omega_\mathrm{b} \) is defined by
 \begin{equation}\label{omega}
 \omega_\mathrm{b} (t)= \omega_0 (1- t/t_\mathrm{crit}) , \mbox{ for } 0\leq t\leq t_\mathrm{crit},
 \end{equation}
 with \(\omega_0= \rho_\mathrm{b} w\)  denoting the initial value of  the blood perfusion,
where \(w\) denotes the volumetric  flow [\si{\per\second}].
The  blood perfusion rate \(\omega_\mathrm{b} \) is set to zero if there is no blood perfusion, which is the case for \(t>t_\mathrm{crit}\).


\section{Analytical solutions}
\label{sanalytical}

Our objective is to build exact solutions that fit, in particular, the physiological problem under study,
\begin{equation}\label{pde}
\tau \partial_t u- \alpha \Delta u + Bu=f,
\end{equation}
via the separable variable method,  in cylindrical coordinates
 \textit{i.e.} at the position \((r,z)\) and the time \(t\).
Using Bernoulli--Fourier technique,  we have 
\begin{equation}\label{odes}
\left(\tau X'(t)+ [ -\alpha(\beta +\eta^2)+ B]X (t)\right)R(r)Z(z)=f(r,z,t),
\end{equation}
where the constants  \(\beta, \eta\in\mathbb{R}\)  are arbitrary,
\begin{description}
\item[Radiative transfer]  \(\tau=1/\nu\), \(\alpha =D\),  \( B=\mu_\mathrm{a}\) and \(f=S\)  in Subsection \ref{s-radiative};

\item[Heat transfer]  \(\tau = \rho c_\mathrm{p}\),  \(\alpha=k\),  \(B= c_\mathrm{b} \omega_\mathrm{b} (t)\) and \(f =q\) in Subsection \ref{s-heat}.
\end{description}
The functions \(X\), \(Z\) and \(R\) are  elementary solutions to the system of ordinary differential equations (ODE)
\begin{equation}\label{ode}
\left\{\begin{array}{l}
X'(t)=\zeta X(t) \\
Z''(z) =\eta^2 Z(z) \\
\left(rR'(r)\right)'=\beta rR(r)  
\end{array} 
\right.
\end{equation}
for some time parameter \(\zeta \in \mathbb{R}\). 
The nonconstant behavior of \(\omega_\mathrm{b} \) (cf. \eqref{omega}) does not invalidate 
the application of this technique in determining an analytical solution for the heat transfer.

Firstly,  a particular solution is available due  to the Duhamel principle (for details, see \cite{lc2024} and the references therein).
Secondly,   let us seek for elementary solutions of the ODE system \eqref{ode}.

The first-order ODE in \eqref{ode}  admits the  elementary solution 
\begin{equation}\label{defX}
X(t) = \exp [ \int_0^t \zeta (s)\dif{s} ], \qquad t>0.
\end{equation}

The second-order ODE in \eqref{ode}, with constant coefficients, admits the elementary solutions
\begin{equation}\label{defZ}
Z(z) = \exp\left[ \pm \eta z \right], \qquad z\in \mathbb{R}.
\end{equation}

The second-order ODE  in \eqref{ode}, with nonconstant coefficients, 
 admits  the Bessel functions of first and second kind and order 0, respectively, \(J_0(\sqrt{|\beta |} r)\) and \(Y_0(\sqrt{|\beta |} r)\)
 if \(\beta <0\);  or the modified Bessel functions of first and second kind and order 0, respectively, 
 \(I_0(\sqrt\beta r)\) and \(K_0(\sqrt\beta r)\) if \(\beta >0\) \cite{ozisik}.
If \(\beta=0\), the elementary solutions reduce to \(R(r) =\log[r]\) and the unity function.

Then, a general solution solving the homogeneous equation \eqref{odes} (\(f=0\))  is available by the above elementary solutions, if
 \begin{equation}\label{zeta}
\tau\zeta +B = \alpha (\beta +\eta^2) .
\end{equation}

Finally, we analyze the PDE \eqref{pde} at the period of time \(0 < t_j+t_p <t<t_{j +1}\), for any \(j = 0,\cdots, N-1\).
This case will describe the homogeneous problem, which is governed without source (\(f=0\)),  at  one pulse-to-pulse interval \(\Delta t\).
For the sake of simplicity, we denote the initial instant of time  \(t_j+t_p\) by \(t_0\) throughout this subsection.

For the initial condition \(u(r,z,t_0)=u_0Z(\eta; z ) \), where 
 \(Z(\eta; z )\) denotes a linear combination of the elementary solutions \(Z\) given as in \eqref{defZ},
we apply the   principle of superposition \cite[Chapter 3-4]{ozisik} to obtain the solution in the form of the Fourier--Bessel series
\begin{equation}\label{useries}
u(r,z,t) = \sum_{m=1}^\infty c_m J_0(b_m r)Z(\eta; z ) \exp[ \zeta_m (t-t_0)]
\end{equation}
that satisfies the homogeneous PDE \eqref{pde}  (\(f=0\)) in a solid cylinder (\(0<r<a\)) with finite length \(L\).
The zero-order Bessel function of the second kind, \(Y_0\),  is excluded  from the solution, because the region
includes the origin \(r = 0\) where \(Y_0\) becomes infinite.

We adapt the homogeneous problems for heat equation (\(B=0\)) in \((r,z,t)\) variables \cite[pp. 127-131]{ozisik},
namely Example 3-9 for a hollow cylinder of finite length and
Example 3-10 for a solid cylinder of  semi-infinite length.

To  determine the coefficients \(c_m\), we multiply  both sides of \eqref{useries}
 by \(rJ_0(b_nr)\) and integrate over the section of the cylinder, obtaining
\[
u_0 \int_0^a r J_0(b_n r) \dif{r} = c_n \int_0^a r J_0^2(b_n r) \dif{r}.
\]
In the left hand side, we used the initial condition 
and in the right hand side, 
we apply the orthogonality relationship of the eigenfunctions 
\begin{equation}
\int_0^a r J_0(b_mr) J_0(b_n r) \dif{r} = \delta_{m,n} \int_0^a r J_0(b_mr) J_0(b_n r) \dif{r}
\end{equation}
where \(\delta_{m,n}\) is the Kronecker delta and the eigenvalues \(b_m\) are determined by the boundary condition.

For any \(a,b>0\), by integration we have \cite[Appendix IV]{ozisik} 
\begin{align*}
\int_0^a rJ_0(br)\dif{r} =& \frac{a}{b}J_1(ba); \\
\int_0^a rJ_0^2(br)\dif{r} =& \frac{a^2}{2}\left(J_0^2 (ba) -J_{-1} (ba) J_1 (ba) \right).
\end{align*}
Recall that  \(W_{-1} (a) = -W_{1} (a) \) for any Bessel function of zero order  \(W= J_0, Y_0, K_0, I_0\).

Then, the complete solution is 
\begin{equation}\label{usol}
u(r,z,t) = \sum_{n=1}^\infty u_0  \frac{2}{b_na} \frac{J_1(b_na)}{J_0^2 (b_na) + J_1^2 (b_na) }J_0(b_n r)Z(\eta; z ) \exp[ \zeta_n (t-t_0)].
\end{equation}


\section{Results}
\label{results}

We focus on the study of the analytical solution  for the fluence rate \(\phi\) in this section.
We derive an exact solution in Section \ref{fluencerate}  and its applicability in Section \ref{profiles}.
This section ends with some considerations for the temperature and respective tissue damage in Section \ref{temperature}.

\subsection{Exact solution for the fluence rate}
\label{fluencerate}

The time parameter in \eqref{zeta} is
\begin{equation}\label{zeta_in}
\zeta  = \nu (D\mu_\mathrm{t} ^2 - \mu_\mathrm{a} ) =
\frac{c}{n} \left(\frac{(\mu_\mathrm{a}+(1-g)^{-1} a(500/\lambda)^b) ^2}{3(\mu_\mathrm{a}+a(500/\lambda)^b)} 
-\mu_\mathrm{a}\right).
\end{equation}
We denote \(\zeta\)  by \(\zeta_\mathrm{in}\)  in the tumor, and \(\zeta_\mathrm{out}\) in the healthy tissue.
According to Table \ref{tabopt},  \(\zeta_\mathrm{in} >0\).

For \(0\leq r\leq r_\mathrm{f}\), \(  0\leq z \leq \ell \) and \( 0\leq t < t_p\),  the transient  solution of \eqref{light}-\eqref{source} is
\begin{equation}\label{phirf}
\phi_f(z,t) \equiv  \frac{ S( r_\mathrm{f} ,z ,0)  }{D\mu_\mathrm{t} ^2 - \mu_\mathrm{a}} \left(  \exp [\zeta_\mathrm{in} t]- 1\right).
\end{equation}
Hereafter, let us set
\begin{equation} \label{defSin}
S_\mathrm{in} =  \frac{ S( r_\mathrm{f} ,0 ,0)  }{D^{(\mathrm{i})}(\mu_\mathrm{t}^{(\mathrm{i})}) ^2 - \mu_\mathrm{a}^{(\mathrm{i})} },
\mbox{  for i = inner}.
\end{equation}

Next,  let us consider the  solution of \eqref{light}-\eqref{source}, 
at two pulses, \(]0; t_p[\) and \(]t_1; t_1 + t_p[\) and the period between them  \(]t_p; t_1[\). 
For the general pulse width \(]t_j; t_j + t_p[\) and the period  \(]t_j +t_p; t_{j+1}[\), we proceed \textit{mutatis mutandis}.

For these two cases, let us extend the fluence rate \(\phi\) to the entire multidomain.
\begin{description}
\item[Case \(0<t<t_p\)] 
The fluence rate \(\phi\), as given in \eqref{phirf}-\eqref{defSin},  verifies the initial condition \(\phi(t=0) =0\).
Moreover, the source \(S\) is a discontinuous function as defined in \eqref{source}-\eqref{planar}.
Indeed, \(S\) can be neglected for \(z>\ell\), as mentioned in Section \ref{profiles}.
Hence, we consistently consider \(\phi=0\) if  \(z>\ell\).

For \(r_\mathrm{f}< r \leq r_\mathrm{o} \) and \(  0\leq z \leq \ell\),  
by the interface continuity condition at \(r=r_\mathrm{f}\) and the boundary condition \eqref{BCz} at \( r_\mathrm{i}\), we may extend \(\phi\) as
\begin{align*}
\phi(r,z,t) &= R_1( r) \exp [ -\mu_\mathrm{t} z]+R_3(r) \exp[ -\mu_\mathrm{t} z + \zeta_\mathrm{in} t ]  \qquad\mathrm{if }\quad r_\mathrm{f}< r \leq r_\mathrm{i};\\
& =  R_2(r)\exp [ -\mu_\mathrm{t} z]+R_4(r) \exp[-\mu_\mathrm{t} z+ \zeta_\mathrm{in} t ]  \qquad\mathrm{otherwise}.
\end{align*}
The existence of \(R_1\) and \(R_3\) rely on  Appendices \ref{particular} and \ref{general}, respectively, 
and \(R_2\) and \(R_4\) are in accordance with Appendix \ref{AppendixB}.
Noting that the fluence rate \(\phi\) verifies the initial condition \(\phi(t=0) =0\), then \(R_1\equiv -R_3\).
This means that
\(S=0\) for \(r>r_\mathrm{f}\) and we consistently might consider \(\phi=0\) either \(r>r_\mathrm{f}\) or \(z>\ell\), if  \( 0\leq t < t_p\).
As we do not expect a zero fluence rate at \(t=t_p \), we assume discontinuous \(\phi\) in time.

The correspondent parameters are \(\beta_1(0)\), \(\beta_1( \zeta_\mathrm{in})\) with \(\beta_1\) being defined in \eqref{beta1},
and  \(\beta_2(0)\), \(\beta_2( \zeta_\mathrm{in})\) with \(\beta_2\) being defined in \eqref{beta2}.
In particular,  we have
\[ 
\beta_2  (\zeta_\mathrm{in})
=\left( \frac{ \mu_\mathrm{a}^{(\mathrm{o})}+a^{(\mathrm{o})}(\lambda/500)^{-b^{(\mathrm{o})}}
}{ \mu_\mathrm{a}^{(\mathrm{i})}+a^{(\mathrm{i})}(\lambda/500)^{-b^{(\mathrm{i})}}}-1\right)(\mu_\mathrm{t}^{(\mathrm{i})})^2 -
 \left( \mu_\mathrm{a}^{(\mathrm{i})}  -  \mu_\mathrm{a}^{(\mathrm{o})} \right)/ D^{(\mathrm{o})},
\] 
where the superscripts (i) and (o) stand for inner and outer regions, respectively.

For \(-\ell < z <0\), the continuity of the fluence rate \(\phi\) and the requirement of \eqref{zeta} being satisfied imply the symmetry of \(\phi\) relative to \(z=0\).

\item[Case \(t_p \leq t < t_1 =  t_p +\Delta t\)]
In this  one pulse-to-pulse interval \(\Delta t\),
 we seek for a solution such that satisfies the homogeneous PDE \eqref{light}
and the initial condition \( \phi_f (z,t_p) \).
Moreover, the solution should decrease in time, with the time parameter \(\zeta\) being defined in \eqref{zeta}.
The principle of superposition now guarantees that the complete solution is constructed as \eqref{usol}, taking
\(\eta =-\mu_\mathrm{t}\) into account. 

For \(0\leq r\leq r_\mathrm{f}\)  and  \(  0\leq z \leq \ell \), we consider  the Fourier--Bessel series
\begin{equation}\label{series}
\phi(r,z,t) = u_0\sum_{n = 1}^\infty  \frac{2}{b_na} \frac{J_1(b_na)}{J_0^2 (b_na) + J_1^2 (b_na) } J_0(\sqrt{|\beta_n|} r) \exp[-\mu_\mathrm{t}z + \zeta_n (t-t_p)].
\end{equation}
The initial constant 
\(u_0= \phi_f(0, t_p) = S_\mathrm{in} \left(  \exp [ \zeta_\mathrm{in} t_p] - 1 \right)\),
 \(a \) denotes the radius correspondent to the under study region of the multidomain,  
and  \(b_n = \sqrt{|\beta_n|}\).
 
For \(0\leq r\leq r_\mathrm{f}\)  and  \(  \ell\leq z \leq L\), 
the function \(\phi\) may be given by
\begin{align*}
\phi(r,z,t) = u_0\sum_{n = 1}^\infty  \frac{2}{b_na} \frac{J_1(b_na)}{J_0^2 (b_na) + J_1^2 (b_na) } J_0(\sqrt{|\beta_n|} r) Z_\ell Z_n(z) \exp[ \zeta_n (t-t_p)], 
\end{align*}
where  
\[
Z_n(z) =\exp[ -\mu_\mathrm{t}^{(\mathrm{i})} \ell ] \frac{\sinh [\eta_n ( L-z ) ]}{ \sinh [\eta_n (L -\ell) ]}.
\]
if the parameters \(\eta_n\) are determined  by \eqref{zeta}, that is, the sequence of parameters verify
\(\eta_n^2 = \left(\zeta_n/\nu + \mu_\mathrm{a}^{(\mathrm{o})}\right)/D^{(\mathrm{o})} - \beta_n >0\),  
for o = outer, or
\[
Z_n(z) =\exp[ -\mu_\mathrm{t}^{(\mathrm{i})} \ell ] \frac{\sin [\eta_n ( L-z ) ]}{ \sin [\eta_n (L -\ell) ]},
\]
if 
\(\eta_n^2 = -\left(  \left(\zeta_n/\nu + \mu_\mathrm{a}^{(\mathrm{o})}\right)/D^{(\mathrm{o})} - \beta_n \right)>0\).

The constant \(Z_\ell\) is determined by  the Robin boundary condition \eqref{BCz} on the function \(\phi\) at \(z =\ell\),  taking  \eqref{defS1} into account.

For the remaining multidomain, we may proceed analogously considering the   principle of superposition \cite{ozisik}.

\end{description}


\subsection{Profiles}
\label{profiles}

In this section, we simulate results by the derived exact solution.
The presented calculations use Octave software, under the optical parameters in Table \ref{tabopt}.

Different wavelengths are used in this work to simulate two  types of lasers: 
Q-switched short pulsed Nd:YAG  laser operating at a wavelength of  \SI{1064}{\nano\metre} 
and  diode short pulsed laser operating at a wavelength of  \SI{810}{\nano\metre}  and  \SI{980}{\nano\metre}.
Also the values of the  pulse widths vary.  
 They can range from \SIrange{500}{2}{\milli\second} in Nd:YAG lasers  \cite{mordon},
 from nanosecond  (\SI{1}{\nano\second} = \SI{1E-9}{\second}) to picosecond (\SI{1}{\pico\second} = \SI{1E-12}{\second})  in most diode lasers
and of femtosecond order (\SI{1}{\femto\second} = \SI{1E-15}{\second}) in free-electron lasers.

\begin{figure}[H]
\begin{multicols}{2}
\centering 
\includegraphics[width=7cm]{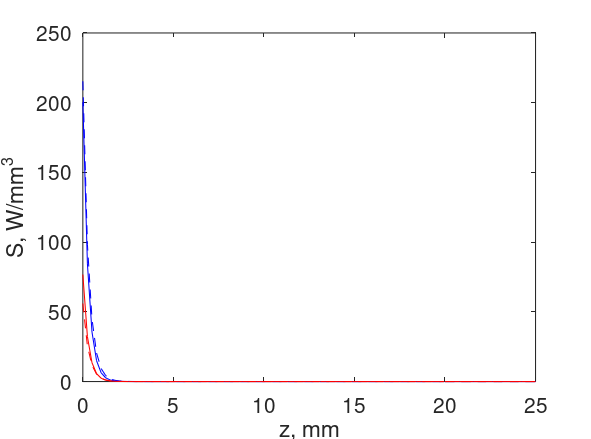} 
\hfill  
\includegraphics[width=7cm]{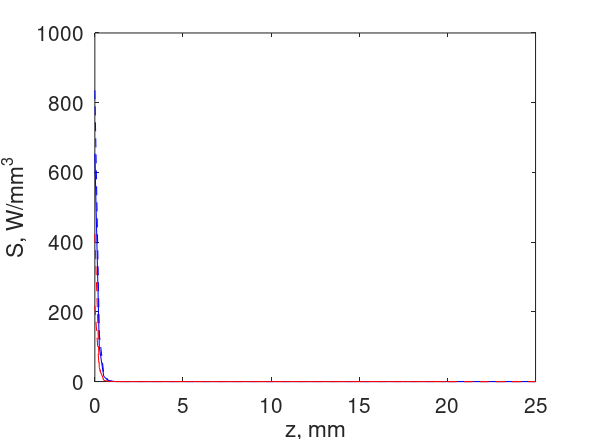}\\
\includegraphics[width=7cm]{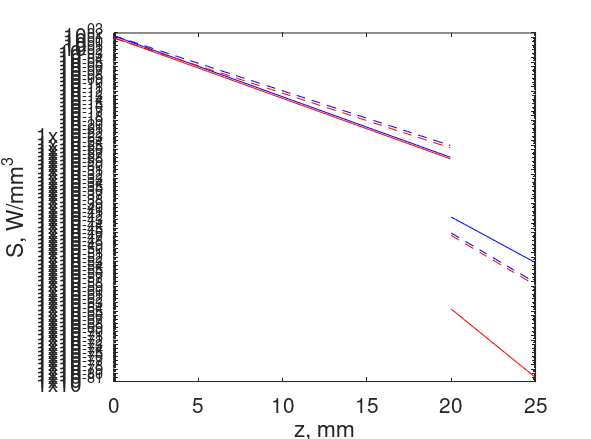}
\hfill  \includegraphics[width=7cm]{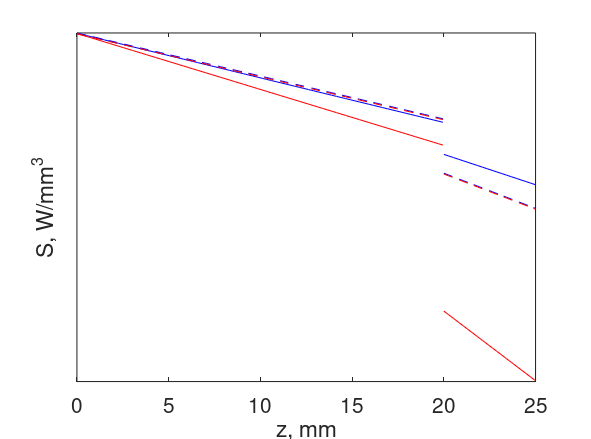}
\end{multicols}
\caption{\footnotesize Graphical representations of the source  \(S\) for the wattage set at \SI{5}{\watt} (in blue):
with wavelengths of  \SI{810}{\nano\metre} (solid line)  and  \SI{980}{\nano\metre} (dashed line)
and at \SI{1.3}{\watt} (in red): with wavelengths of   \SI{980}{\nano\metre} (dashed line) and  \SI{1064}{\nano\metre} (solid line).
 (\textbf{a})   Plot with linear axes of the breast source.  (\textbf{b})  Plot using a logarithmic scale for the \(y\)-axis of the breast source.  
 (\textbf{c}) Linear plot of the prostate source.  (\textbf{d}) Semilog plot of the prostate source. For details see Table~\ref{tabmax}.
\label{sourceS}} 
\end{figure}

We begin by analyzing the light absorption of the tissues according to the BL law.
  Figure~\ref{sourceS} illustrates the linear and semi-log steady-state profiles of \(S\) for these operating wavelengths, 
under optical values of the carcinoma-adipose breast and tumor-healthy prostate tissues
 according to Table~\ref{tabopt}.

At the tip \(z=0\), the source \(S(0, 0, t_p) \) has its maximum values according to Table~\ref{tabmax}.
\begin{table}[H]
\caption{Source values \(S\)  at \(z=0\) (maximum), and at the discontinuity \(z=\ell\) [\si{\watt\per\cubic\milli\metre}].
\label{tabmax} } 
\begin{tabular}{lcccc}
\toprule
\(\lambda\)  [\si{\nano\metre}]  & 810  & 980 & 980 &1,064 \\ 
\(P_\mathrm{peak}\)  [\si{\watt}] & 5& 5 & 1.3 & 1.3 \\
\midrule
Maximum breast tumor & 198   & 214 & 56   & 77  \\
Minimum breast tumor & 9.5e-28    &  6.8e-25  & 1.8e-25   & 3.9e-28  \\
Maximum adipose breast &4.4e-42    & 6.3e-46   & 1.6e-46 & 3.1e-64 \\
Maximum prostate tumor & 647 & 828 & 215 & 420 \\
Minimum prostate tumor &  6.6e-62 & 2.2e-59 & 5.7e-60  & 2.7e-78 \\
Maximum healthy prostate &2.6e-86    & 6.6e-100   & 1.7e-100 & 6.8e-200 \\
\bottomrule
\end{tabular}
\end{table}

We emphasize that the source \(S\) does not have expression in the outer region  correspondent to the healthy tissue.
Then, we may neglect the source  \(S\) for \(z>\ell\). 
For instance, 
at  the wattage set at \SI{5}{\watt} (in blue), with wavelength of  \SI{810}{\nano\metre} (solid line),
we have   \(S(0, 20, t_p)\simeq \SI{e-27}{\watt\per\cubic\milli\metre}\) for the breast tumor border and \(S(0, 20, t_p)\simeq \SI{e-42}{\watt\per\cubic\milli\metre}\) 
for the adipose breast border.

\begin{figure}[H]
\begin{multicols}{2}
\centering 
\includegraphics[width=7cm]{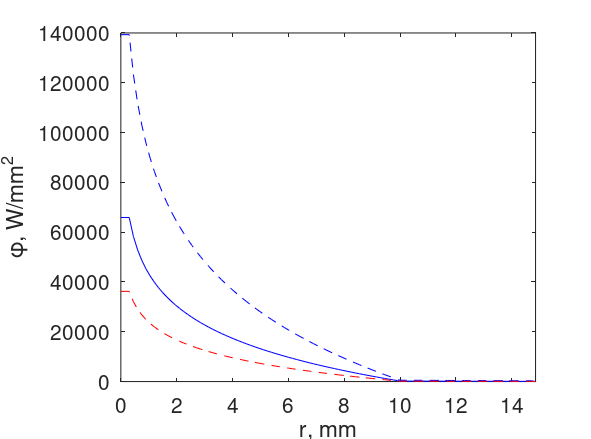} \hfill  \includegraphics[width=7cm]{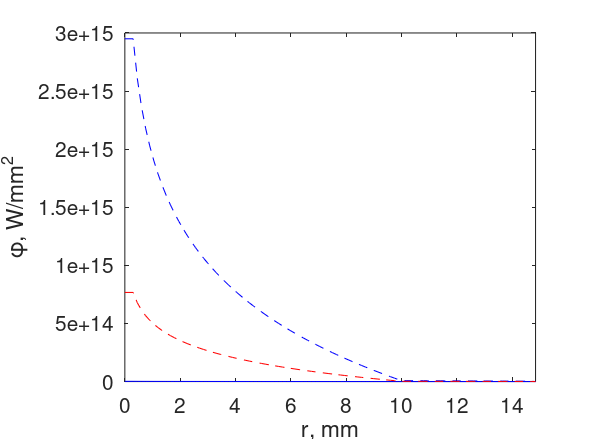}\\
\includegraphics[width=7cm]{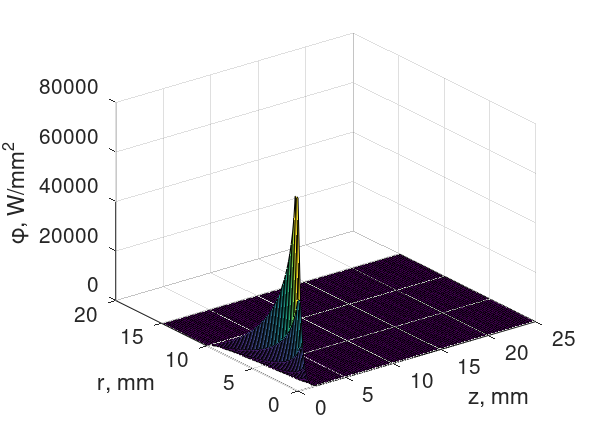} \hfill  \includegraphics[width=7cm]{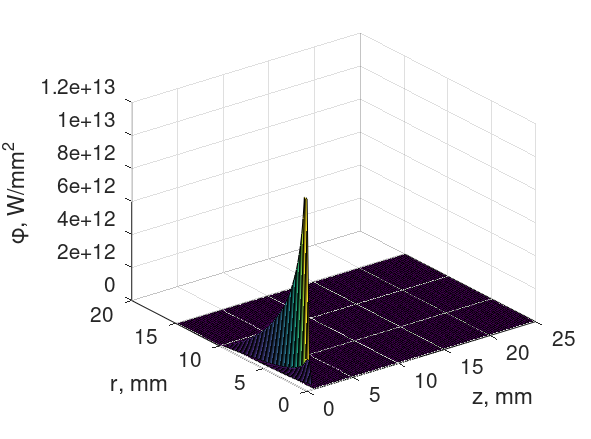}
\end{multicols}
\caption{\footnotesize (\textbf{a})  Breast radial graphical representations of the fluence rate  \(\phi\) for the
wavelength of  \SI{810}{\nano\metre} and wattage set at \SI{5}{\watt} (blue solid line)  and 
for  \SI{980}{\nano\metre} and  wattage set at \SI{5}{\watt} (blue dashed line)
and at \SI{1.3}{\watt} (red dashed line).
 (\textbf{b})   Plot for the tumor-adipose breast tissue for the
wavelength of  \SI{980}{\nano\metre} and wattage set at \SI{1.3}{\watt}.
 (\textbf{c}) Prostate  radial graphical representations of the fluence rate  \(\phi\) as in (a).
 (\textbf{d})  Plot for the tumor-healthy prostate tissue for the
wavelength of  \SI{810}{\nano\metre} and wattage set at \SI{5}{\watt}.
\label{radial}} 
\end{figure}

The fluence rate \(\phi\) has different graphical representations displayed in Figure~\ref{radial},  at diode pulse width \(t_p=\SI{ 1e-11}{\second}\).
Figure~\ref{radial} shows distributions of \(\phi\) for breast in (a) and prostate in (c),
as function of the radius \(r\),  at the different operating  systems.
Figure~\ref{radial} (b)  shows distributions of \(\phi\) for breast for
 the wattage set at \SI{1.3}{\watt} and  at a wavelength of  \SI{980}{\nano\metre} corresponding to the red dashed line in  (a).
Figure~\ref{radial} (d)  shows distributions of \(\phi\) for prostate for
 the wattage set at \SI{5}{\watt} and  at a wavelength of  \SI{810}{\nano\metre} corresponding to the the blue solid line in  (b).
At the referred values of pulse width in the order of \si{\milli\second} for Nd:YAG, the fluence rate \(\phi\) 
does not fit a coherent result.  Indeed, the value of \(\exp[\zeta_\mathrm{in} t_p] \gg 1e+304\).
This drawback is due to the presence of the   light velocity \(\nu\) in the tissue.

\begin{figure}[H]
\begin{multicols}{2}
\centering 
\includegraphics[width=7cm]{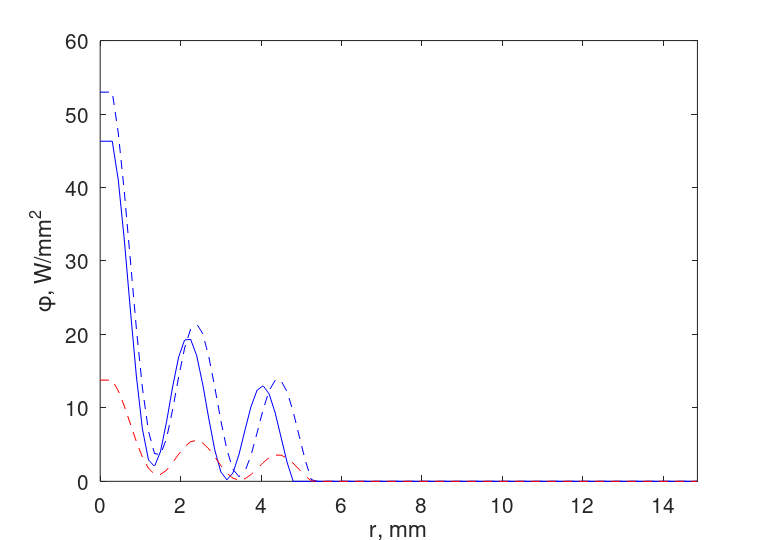}  \hfill    \includegraphics[width=7cm]{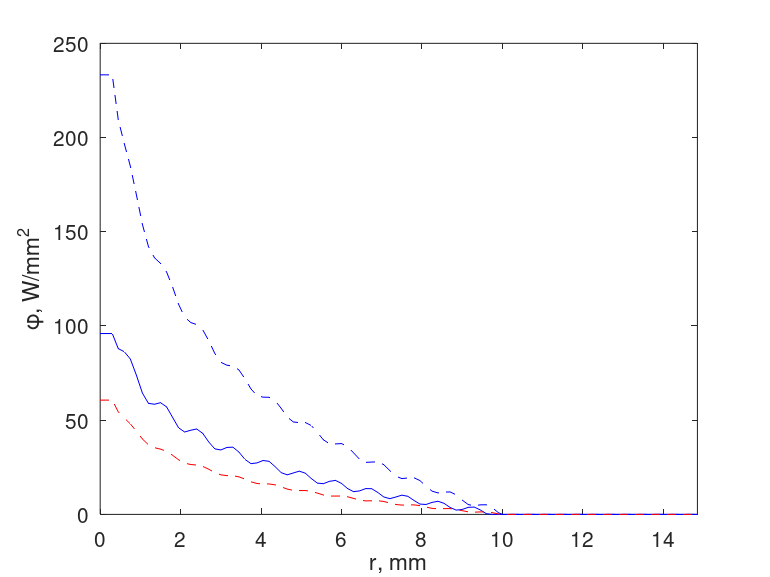}\\
\includegraphics[width=7cm]{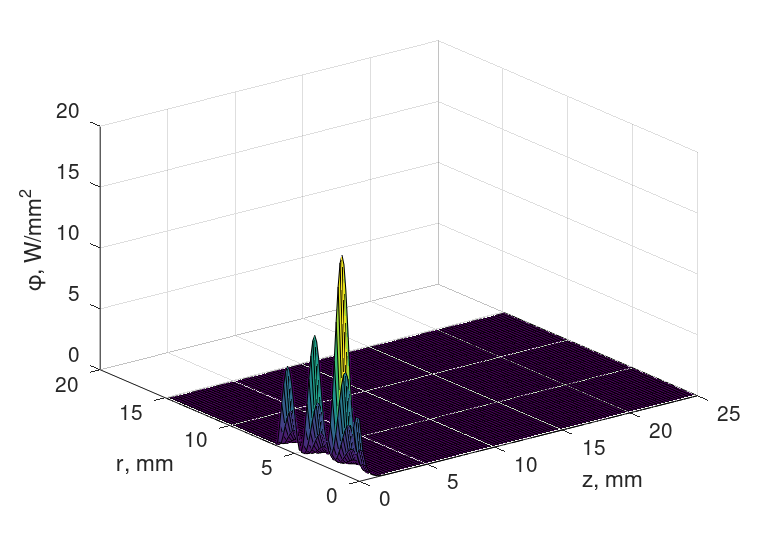}  \hfill  \includegraphics[width=7cm]{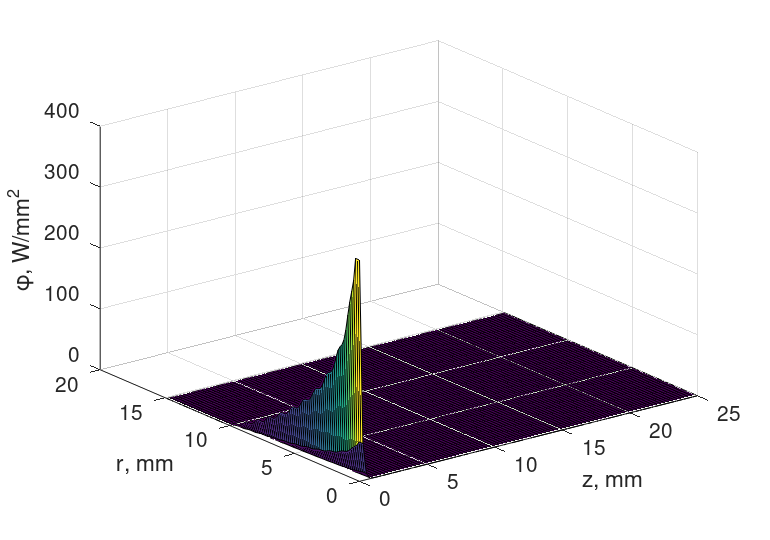}
\end{multicols}
\caption{\footnotesize (\textbf{a})  Breast radial graphical representations of the fluence rate  \(\phi\) for the
wavelength of  \SI{810}{\nano\metre} and wattage set at \SI{5}{\watt} (blue solid line)  and 
for  \SI{980}{\nano\metre} and  wattage set at \SI{5}{\watt} (blue dashed line)
and at \SI{1.3}{\watt} (red dashed line).
 (\textbf{b})   Plot for the tumor-adipose breast tissue for the
wavelength of  \SI{980}{\nano\metre} and wattage set at \SI{1.3}{\watt}.
 (\textbf{c}) Prostate  radial graphical representations of the fluence rate  \(\phi\) as in (a).
 (\textbf{d})  Plot for the tumor-healthy prostate tissue for the
wavelength of  \SI{980}{\nano\metre} and wattage set at \SI{1.3}{\watt}.
\label{rr3D}} 
\end{figure}

Next,  Figure~\ref{rr3D} displays  graphical representations of the fluence rate \(\phi\),  at diode pulse width \(t_p=\SI{ 1.0}{\pico\second}\).
Figure~\ref{rr3D} shows distributions of \(\phi\) for breast in (a) and prostate in (c),
as function of the radius \(r\),  at the different operating  systems.
Again, the result relative to the situation of wavelength of  \SI{1,064}{\nano\metre}  is excluded by the same reason as before.
Figure~\ref{rr3D} (b) and (d) show distributions of \(\phi\) for breast and  prostate, respectively,  for
 the wattage set at \SI{1.3}{\watt} and  at a wavelength of  \SI{980}{\nano\metre} corresponding to the red dashed line in  (a) and  (b).
At  the temporal pulse width \(t_p\), 
the fluence rate does not  rapidly  decrease with distance from the tip surface, which implies that the
 the Joule effect heating is generated along the \(r\)-direction.
Also in Figure~\ref{rr3D} (b) and (d),
the slopes of curves in the \(z\)-direction illustrate that  the fluence rate \(\phi\) is locally concentrated on the origin  relative to the tip surface.
Distributions of the fluence rate over time until \(t_\mathrm{end}\) are subject of study in forthcoming work.
 The analysis of the behavior of the outer radius \(r_\mathrm{o}\)  is still an open problem.

\subsection{Temperature and tissue damage}
\label{temperature}

The  time dependent  parameter \(\zeta\) in \eqref{zeta} is
\[ \zeta(\beta; t) = 
\left(  \rho c_\mathrm{p}\right)^{-1} \left( k(\beta +\mu_\mathrm{t} ^2 ) - c_\mathrm{b} \rho_\mathrm{b}w(1- t/t_\mathrm{crit})\right).
\]
We set 
\[
\zeta_1(\beta; t) = \left\{ \begin{array}{ll}
\zeta(\beta;0) &\mbox{ for }t  < t_p\\
\zeta (\beta; t) &\mbox{ for } t_p \leq t < t_\mathrm{crit}\\
\zeta(\beta;t_\mathrm{crit}) &\mbox{ for }t  \geq t_\mathrm{crit}
\end{array} \right.
\] 
under the parameter \(\beta\), namely \(\beta =0\),  \(\beta =\beta_1\)  or \(\beta =\beta_2\),  according to Appendices \ref{AppendixA}  and \ref{AppendixB}.
Notice that \(\zeta\) increases as function of time and \(\zeta_1\) obeys \( \zeta(\beta;0) \leq \zeta_1(\beta, t)  \leq\zeta(\beta; t_\mathrm{crit})\), for all \(t>0\).
For the sake of simplicity,  we write the constant \(\zeta_0 = \zeta(0; 0)\).
According to Table~\ref{tabthem},  we can calculate \(\zeta_0\) in the tumors (cf. Table~\ref{tabz0}).
\begin{table}[H]
\caption{Average thermal parameters for blood and for breast and prostate (tumor and healthy) tissues \cite{ho,marq}. \label{tabthem}     }  
\begin{tabular} {llccccc}
\toprule
 & \textbf{unit} & \textbf{blood} & \textbf{breast tumor} & \textbf{prostate tumor} &\textbf{gland} & \textbf{fatty tissue} \\
 \midrule
  \(\rho\) &\si{\kilogram\per\cubic\metre}  &    1060.00  &1000.00 &999.00 & 1041.00 & 920.00  \\
\(\omega_0\)  & \si{\kilogram\per\cubic\metre\per\second}  & & 0.5 & 0.416& 0.54 & 0.32   \\
\bottomrule
\end{tabular}
\end{table}
The breast tumor is commonly located in the gland tissue.
If the breast tumor is located in the adipose tissue,  the healthy tissue will exhibit
 lower  density and  blood perfusion values than ones of the tumor, as for the prostate values (cf. Table~\ref{tabthem}).
Also, the density of any tumor varies from lower values at a primary stage and then increases at very high values.
\begin{table}[H]
\caption{Values for parameter \(\zeta_0\) at different wavelengths. \label{tabz0}}       
\begin{tabular}{lcc}
\toprule
\(\lambda\)  [\si{\nano\metre}] & \textbf{breast}&\textbf{prostate} \\
\midrule
810 & 1.7e+3 & 7.9e+3\\
980&   1.4e+3 & 7.3e+3 \\
1,064&   1.7e+3& 1.2e+4 \\
\bottomrule
\end{tabular}
\end{table}

Thanks to the Duhamel principle (for details see \cite{lc2024}),
the  solution of \eqref{bioheat}-\eqref{power} is, for  \(  0\leq r\leq r_\mathrm{f}\) and \( 0\leq z \leq \ell \),
\begin{itemize}
\item if  \( 0\leq t < t_p\):
\begin{align*}
T(r,z,t) &= T_\mathrm{b} +  \frac{ \mu_\mathrm{a}  }{\rho c_\mathrm{p}}
\frac{ S( r_\mathrm{f} ,z ,0)  }{D\mu_\mathrm{t} ^2 - \mu_\mathrm{a}}    
\int_0^t\left(  \exp [\zeta_\mathrm{in} s]- 1\right)\exp[\zeta_0 (t-s)]\dif{s}\\
 &=  T_\mathrm{b} +  \mu_\mathrm{a}  S_\mathrm{in} \exp [- \mu_\mathrm{t} z] \left(
\frac{\exp [\zeta_\mathrm{in} t] - \exp [\zeta_0 t]}{\rho c_\mathrm{p} (\zeta_\mathrm{in} - \zeta_0)} 
+ \frac{1-\exp [\zeta_0 t]}{ k \mu_\mathrm{t} ^2- c_\mathrm{b}\rho_\mathrm{b} w} \right) ;
\end{align*}
\item  if \( t \geq t_p\):
\[
T(r,z,t) = T(r_\mathrm{f},z,t_p) +  \frac{ \mu_\mathrm{a}  }{ c_\mathrm{p}\rho}
\int_{t_p}^t \phi( r_\mathrm{f} ,z ,s) \exp [\int_{0}^{t-s}\zeta_1 (0; \tau)\dif{\tau}]\dif{s}.\]
\end{itemize}
We may extend the temperature \(T\) to the entire multidomain, by proceeding as in Section \ref{fluencerate}.

To evaluate the thermal damage, we need to estimate \(t_\mathrm{crit}\). From \eqref{arr},  we have
\begin{equation}\label{crit}
\frac{1}{A}\exp\left[\frac{E_\mathrm{a} }{ RT_\mathrm{max}}\right] \leq t_\mathrm{crit}\leq \frac{1}{A} \exp\left[\frac{E_\mathrm{a} }{ RT_\mathrm{min}}\right],
\end{equation}
where
\[ 
T_\mathrm{max} =\sup_{  \Omega\times ]0,t_\mathrm{crit}[}  T;\qquad
T_\mathrm{min} =\inf_{ \Omega\times ]0,t_\mathrm{crit}[}  T.
\]

The minimum  temperature \(T_\mathrm{min}=T_\mathrm{b}\) and  \(t_\mathrm{crit}\) obeying \eqref{crit} verifies
\begin{equation}\label{tmax}
 t_\mathrm{crit}\leq \frac{1}{A} \exp\left[\frac{E_\mathrm{a} }{ RT_\mathrm{b}}\right].
\end{equation}

 Typical values for prostate tumor are \(E_\mathrm{a} = \SI{5.67E5 } {\joule\per\mol}  \) and 
\(A = \SI{1.7E91}{\per\second}  \) \cite{HeBischof}.
Based on these values, \eqref{tmax} infers that  \(t_\mathrm{crit}\) obeys the upper bound \( t_\mathrm{crit}\leq \SI{9.9E+03}{\second}\sim \SI{2.75}{\hour}\).
As commonly accepted, the maximum temperature  with respect to produce complete in situ destruction of tissue 
is  \(T_\mathrm{max}=\SI{50}{\degreeCelsius} \).
Then,   \(t_\mathrm{crit}\) similarly obeys  the lower bound   \( t_\mathrm{crit}\geq \SI{2.8871}{\second}\).


\section{Discussion and conclusions}
\label{discuss}

We firstly  emphasize that the validness  of the diffusion  approximation of the radiative transfer equation \eqref{light} is true.
Some important findings include the fact that
the criterion \(\mu_\mathrm{a} \ll \mu'_\mathrm{s}\) can not be disregard in the validity of the stationary diffusion approximation 
as stated but not proven  in  \cite{w-vang6} for any geometric domain: plane \cite[Section 6.4.1]{w-vang6},  spherical with isotropic point source \cite[Section 6.4.2]{w-vang6}, 
or cylindrical  \cite[Section 6.4.3]{w-vang6}. 
The criterion \(\mu_\mathrm{a} \ll \mu'_\mathrm{s}\)  in the validity is precisely
quantified in \cite{Capart} as \(\mu_\mathrm{a} /\mu'_\mathrm{s}\sim 1/100\).  
According to Table \ref{tabopt}, we calculate \(\mu_\mathrm{a}/ \mu'_\mathrm{s}\) (cf. Table~\ref{tabv}).
\begin{table}[H]
\caption{Values for  \(\mu_\mathrm{a}/ \mu'_\mathrm{s}\) at different wavelengths. \label{tabv}}       
\begin{tabular}{*{27}{c}} 
\toprule
\(\lambda\)  [\si{\nano\metre}] & \textbf{breast tumor} & \textbf{breast tissue}&\textbf{prostate tumor}&\textbf{prostate tissue} \\
\midrule
810 &7.9e-03  & 1.7e-02 &  8.2e-03 &  4.2e-02\\
980 &   9.2e-03 & 2.4e-02 &  1.0e-02   & 4.7e-02 \\
1,064 &  8.9e-03 & 4.0e-02 &  1.1e-02 &  4.3e-02 \\
\bottomrule
\end{tabular}
\end{table}
The criterion is verified by  all values. Although some of the values 
are indeed greater, they verify  \(\mu_\mathrm{a} /\mu'_\mathrm{s}\sim 1/100\).  
In contrast in \cite{CarpPrahl},
the authors argue that the \(\delta\)-Eddington approximation to the Boltzmann transport equation provides small percentage error if 
\(\mu'_\mathrm{s}/\mu_\mathrm{a} \sim 100\).

Our first finding is that the source \(S\) does not reach the outer region  correspondent to the healthy tissue,
when the laser emission tip surface is located at the middle of the tumor, for one pulse width.
  Figure~\ref{sourceS} and Table~\ref{tabmax}  confirm the local behavior of the source \(S\).
 Its longitudinal discontinuity shown  in   Figure~\ref{sourceS} (b) and (d) is a consequence of the tissue dependence on the optical values. 
Therefore, if the laser tip is fix at the center of the tumor, its action is local.
Then, a moving tip, as used in EVLA \cite{lc2024}, should be object of study in future research in order to find out
if  a better performance will be provided.

The relationship \eqref{zeta} implies some others conclusions about  the diffusion approximation of the radiative transfer equation \eqref{light}.
If source of scattered photons  \(S\) is given in function of the   effective attenuation coefficient
 \(\mu_\mathrm{eff} = \sqrt{ 3 \mu_\mathrm{a} (\mu_\mathrm{a}+\mu'_\mathrm{s}) }\) \cite{oden}, the  relationship \eqref{zeta} reads 
\( \zeta/\nu + \mu_\mathrm{a}=  D\mu_\mathrm{eff}^2 = \mu_\mathrm{a}\). This means 
 that the fluence rate will be steady-state, because \(\zeta =0\). 
In \cite{w-vang7}, the author presents a detailed analytical solution of the  diffuse radiant fluence rate for time-independent searchlight problem
 on flat-slab and semi-infinite geometries, under homogeneous Dirichlet conditions, by using the Green's function method. 
Also, the author refers some inaccuracies of the diffusion approximation.

In Section \ref{fluencerate}, we derive the fluence rate \(\phi\) as a solution of the  diffusion approximation of the radiative transfer equation \eqref{light}.
 Its  simulations, throughout the inlet region, originated by the tumor,  and  the outlet region, originated by the healthy tissue,
are illustrated in Figures \ref{radial} and \ref{rr3D},  having different pulse widths.

On the one hand,
we observe different behaviors occur at different pulse width values by comparing Figures \ref{radial} and \ref{rr3D}, 
which suggest that there are values of the pulse width that blow up the exponential  function, namely
at  pulse widths \(t_p=  \SI{ 10}{\pico\second}\) or  above.
However,  at pulse widths of order   \si{\pico\second} we may conclude that a
further study  of the  time-dependent  searchlight problem is a priority.

On the other hand, 
the present results in Figure~\ref{rr3D} (a) and (c) suggest that power alone  is a stronger contributor of the profiles magnitude
than wavelengths.
We observe that Figure~\ref{radial} (c) for the prostate tissue indicates that it is the wavelength the main contributor  of the profiles magnitude.
However, taking into account the previous comment about the comparison between  Figures \ref{radial} and \ref{rr3D},
  we discard the last observation in order to avoid ambiguities and erroneous conclusions.

Our last conclusion is that the treatment should  have a moving focal point with short exposure time \(t_\mathrm{end} =t_\mathrm{crit}\)  to preserve the healthy tissue.
We emphasize that to evaluate \(t_\mathrm{crit}\), as defined in Section \ref{temperature} by  the Arrhenius burn integration \eqref{arr},
is based on the dependence of the parameters  \(A\) or \(E_a\).
A breaking point for either \(A\) or \(E_a\) is reported in \cite{ryl}, at \(T=\SI{54}{\degreeCelsius} \), which is greater than \SI{50}{\degreeCelsius}
(known as the maximum temperature  with respect to produce complete in situ destruction of tissue).
To produce the complete cure of the patient at a least treatment time (for instance \(t_\mathrm{crit} = \SI{2.8871}{\second}\), for prostate tumor) is the desired objective.

To the study of the heat exchange problem for ablation performed in biomedical sciences,
it is important that its source, the absorbed optical power density \eqref{power},   be analytically precised.
The various equations in the present paper allow to calculate the light fluence rate as a function of \((r,z,t)\), within different regions in space. 
We apply them to breast and prostate tumors, but they may be applied to any kind,
by recurring to the optical and thermal properties of the tissue in study. 
According to the experimental methods,  we set higher optical values for the tumor tissue than ones of the healthy, 
but our results may be replicated to inverse optical values.
 

\subsection*{Data availability}
The data supporting presented results is publicly available in the repository: LuisaConsiglieri/FocalLaserAblationStudy (2025).

\href{https://doi.org/10.5281/zenodo.14993810}{https://doi.org/10.5281/zenodo.14993810}

\subsection*{Acknowledgements}
Grateful thanks to Daniel Pfaller for his  comments on the laser modeling.


\appendix
\numberwithin{equation}{section}

\section[\appendixname~\thesection]{Extending inside the tumor}
\label{AppendixA}

Let us consider 
\[
\phi(r,z,t) =  R( r_\mathrm{f} ) \exp [  -\mu_\mathrm{t}  z+  \zeta t] ,
\]
where \(R( r_\mathrm{f} )\) is a positive constant, which stands for
\begin{itemize} 
\item  the particular solution of  \eqref{light}-\eqref{source}, if \(\zeta = 0\);
\item the general solution of  \eqref{light}-\eqref{source}, if \( \zeta  = \nu (D\mu_\mathrm{t} ^2 - \mu_\mathrm{a} ) \) (or \( \zeta = \zeta_\mathrm{in}\), cf. \eqref{zeta_in}).
\end{itemize}

For  \(r_\mathrm{f}< r \leq r_\mathrm{i}\),  we seek a function in the form
\begin{equation}\label{phio}
\phi(r,z,t) = R_1(r)\exp [-\mu_\mathrm{t}  z +  \zeta t],
\end{equation}
with \(R_1\) being a solution to the second-order ODE  in \eqref{ode}, with nonconstant coefficients, with parameter \(\beta_1\) verifying
\begin{equation}\label{beta1}
\beta_1 (\zeta ) = \left(\zeta/\nu + \mu_\mathrm{a}^{(\mathrm{i})}\right)/D^{(\mathrm{i})} - (\mu_\mathrm{t}^{(\mathrm{i})})^2 \quad \mbox{ for i = inner}.
\end{equation}

\subsection[\appendixname~\thesubsection]{Particular solution}
\label{particular}

In this Section, \(\beta_1\) reduces to
\[
\beta_1 = \mu_\mathrm{a}^{(\mathrm{i})}/D^{(\mathrm{i})} - (\mu_\mathrm{t}^{(\mathrm{i})})^2 <0.
\]
According to Table~\ref{tabopt}, the constant \(\beta_1\) is negative.

Using the continuity condition  on the function \(\phi\) at \(r=r_\mathrm{f}\), we have
\begin{equation}
R_1(r) = C_{1} J_0(\sqrt{|\beta _1|} r) +
 \frac{ R( r_\mathrm{f} ) - C_{1} J_0(\sqrt{|\beta _1|} r_\mathrm{f}) }{Y_0(\sqrt{|\beta _1|} r_\mathrm{f}) }Y_0(\sqrt{|\beta _1|} r) ,
\mbox{ for }C_1\in\mathbb{R}.\label{defR1}
\end{equation}

The  involved constant  \(C_1\) in \eqref{defR1} may be determined by  the continuity condition on the flux
\begin{equation}\label{C1}
C_1  J_1(\sqrt{|\beta _1|} r_\mathrm{f}) +
 \frac{R( r_\mathrm{f})  - C_{1} J_0(\sqrt{|\beta _1|} r_\mathrm{f}) }{Y_0(\sqrt{|\beta _1|} r_\mathrm{f}) }Y_1(\sqrt{|\beta _1|} r_\mathrm{f})=0.
\end{equation}

Applying the wronskian relationship   (see \cite[p. 672]{ozisik} and \cite[pages 360 and 375]{olver})
\begin{equation}\label{ozisik}
J_1(\beta r)Y_0(\beta r)-Y_1(\beta r)J_0(\beta r) =\frac{2}{\pi\beta r},
\end{equation}
for  \(\beta = \sqrt{|\beta _1|}>0\), we may compute \eqref{C1}, concluding 
\[ 
C_1 = -  \frac{\pi  }{2 } \sqrt{|\beta _1|} r_\mathrm{f} R( r_\mathrm{f}) Y_1(\sqrt{|\beta _1|} r_\mathrm{f}) .
\] 

\subsection[\appendixname~\thesubsection]{General solution}
\label{general}
 
Notice that the relationship \eqref{zeta_in} between the parameters \(\zeta\) and \(\mu_\mathrm{t}\) leads \(\beta_1 =0\), then 
\(R_1(r) = R(r_\mathrm{f}) - \beta_0\log[ r/r_\mathrm{f} ]\), which verifies  \(R_1(r_\mathrm{f})= R(r_\mathrm{f})\).
Here, \(b_0\) is a positive constant determined by the Robin boundary condition \eqref{BCr}, namely it is given by
\[
b_0 = \frac{\gamma_r}{ \gamma_r \log[ r_\mathrm{i} /r_\mathrm{f}] -2D/r_\mathrm{i} } R(r_\mathrm{f}).\]

\section[\appendixname~\thesection]{Extending outside the tumor (\(z<\ell\))}
\label{AppendixB}

In the tumor, let us consider the general solution
\[
\phi(r,z,t) =  R_1(r) \exp [  -\mu_\mathrm{t}  z+\zeta t] ,
\]
with \(R_1\) being a solution to the second-order ODE  in \eqref{ode}, with nonconstant coefficients,, with parameter \(\beta_1 = \beta_1 (\zeta ) \)
 verifying \eqref{beta1}.

For \(r_\mathrm{i}< r \leq r_\mathrm{o}\), we seek a function in the form
\[ 
\phi(r,z,t) =   R_2(r) \exp [  -\mu_\mathrm{t}  z+ \zeta t] ,
\] 
with \(R_2\) being constituted by Bessel functions with parameter \(\beta_2\) such that
\[
\beta_2= \left(\zeta/\nu + \mu_\mathrm{a}^{(\mathrm{o})}\right)/D^{(\mathrm{o})} - (\mu_\mathrm{t}^{(\mathrm{i})})^2 \quad \mbox{ for o = outer}.
\]
Considering \eqref{beta1}, \(\beta_2\) reads
\begin{equation}\label{beta2}
\beta_2 (\zeta) =\left\{ \begin{array}{ll} 
\mu_\mathrm{a}^{(\mathrm{o})} /D^{(\mathrm{o})} - (\mu_\mathrm{t}^{(\mathrm{i})})^2   &\mbox{if } \zeta = 0\\
\left( D^{(\mathrm{i})} \beta_1 +  \left( D^{(\mathrm{i})} - D^{(\mathrm{o})}   \right) (\mu_\mathrm{t}^{(\mathrm{i})})^2 -
\mu_\mathrm{a}^{(\mathrm{i})}  +  \mu_\mathrm{a}^{(\mathrm{o})} \right)  / D^{(\mathrm{o})} & \mbox{ if } \zeta = \zeta_\mathrm{in}.
\end{array}\right.
\end{equation}

\subsection[\appendixname~\thesubsection]{Case \(\beta_2 <0\) (Bessel functions)}
\label{AppendixB1}

Using the continuity condition  on the function \(\phi\)  at \(r=r_\mathrm{i}\), we have
\[ 
R_2(r) = C_{2} J_0(\sqrt{|\beta _2|}  r) +  \frac{ R_1(r_\mathrm{i})  - C_{2} J_0(\sqrt{|\beta _2|}  r_\mathrm{i}) }{Y_0(\sqrt{|\beta _2|}  r_\mathrm{i}) }Y_0(\sqrt{|\beta _2|}  r),
\quad C_2\in\mathbb{R} .
\] 
The constant  \(C_2\) may be determined by  the  Robin boundary condition \eqref{BCr}, namely
\begin{equation}\label{C2}
C_{2} J_1(\sqrt{|\beta _2|}  r_\mathrm{i}) + 
 \frac{  R_1(r_\mathrm{i}) - C_{2} J_0(\sqrt{|\beta _2|}  r_\mathrm{i}) }{Y_0(\sqrt{|\beta _2|}  r_\mathrm{i}) }Y_1(\sqrt{|\beta _2|}  r_\mathrm{i})
= - \frac{1}{\sqrt{|\beta_2|}}  \frac{D^{(\mathrm{i})} }{D^{(\mathrm{o})} } R_1'(r_\mathrm{i})  .
\end{equation}
Taking the wronskian relationship  \eqref{ozisik} into account, 
for  \(\beta = \sqrt{|\beta _2|}>0\), we may compute \eqref{C2}, concluding 
\begin{equation}\label{defC2}
C_{2} =  - \frac{\pi  }{2 }  r_\mathrm{i}\left(   \frac{D^{(\mathrm{i})} }{D^{(\mathrm{o})} } R_1'(r_\mathrm{i}) Y_0(\sqrt{|\beta _2|}  r_\mathrm{i}) +
\sqrt{|\beta _2|}  R_1(r_\mathrm{i}) Y_1(\sqrt{|\beta _2|}  r_\mathrm{i})  \right) .
\end{equation}

\subsection[\appendixname~\thesubsection]{Case \(\beta_2 >0\) (modified Bessel functions)}
\label{AppendixB2}

Using the continuity condition  on the function \(\phi\)  at \(r=r_\mathrm{i}\), we have
\[
R_2(r) = C_{2} I_0(\sqrt{\beta _2} r) +
 \frac{ R_1( r_\mathrm{i} ) - C_{2} I_0(\sqrt{\beta _2} r_\mathrm{i}) }{K_0(\sqrt{\beta _2} r_\mathrm{i}) }K_0(\sqrt{\beta _2} r) ,
\quad C_2\in\mathbb{R} .\] 

The wronskian relationship for modified Bessel functions (see  \cite[pages 360 and 375]{olver}) is 
\begin{equation}
 I_1(\beta r) K_0(\beta r) + K_1(\beta r) I_0(\beta r) =\frac{1}{\beta r} .\label{wronsk2}
\end{equation}
Applying \eqref{wronsk2} for \(\beta = \sqrt{\beta _2} >0\), analogously to \eqref{C2}-\eqref{defC2},  we conclude 
\[
C_{2} =    r_\mathrm{i}\left(  \frac{D^{(\mathrm{i})} }{D^{(\mathrm{o})} }
 R_1'(r_\mathrm{i}) K_0(\sqrt{\beta _2}  r_\mathrm{i}) + \sqrt{\beta _2}  R_1(r_\mathrm{i}) K_1(\sqrt{\beta _2}  r_\mathrm{i}) \right) .
\]



\bibliographystyle{plain}

 \bibliography{FLAbiblio_abbrv.bib}

\begin{thebibliography}{10}

\bibitem{bane}
A.~Banerjee, A.A. Ogale, C.~Das, K.~Mitra, and C.~Subramanian.
\newblock Temperature distribution in different materials due to short pulse
  laser irradiation.
\newblock {\em Heat Transfer Eng.}, 26(8):41--49, 2005.

\bibitem{merryman}
S.M. Boronyak and W.D. Merryman.
\newblock Development of a simultaneous cryo-anchoring and radiofrequency
  ablation catheter for percutaneous treatment of mitral valve prolapse.
\newblock {\em Ann. Biomed. Eng.}, 40(9):1971--1981, 2012.

\bibitem{Capart}
A.~Capart, S.~Ikegaya, E.~Okada, M.~Machida, and Y.~Hoshi.
\newblock Experimental tests of indicators for the degree of validness of the
  diffusion approximation.
\newblock {\em J. Phys. Commun.}, 5(025012), 2021.

\bibitem{CarpPrahl}
S.A. Carp, S.A. Prahl, and V.~Venugopalan.
\newblock {Radiative transport in the delta-\(P_1\) approximation: Accuracy of
  fluence rate and optical penetration depth predictions in turbid
  semi-infinite media}.
\newblock {\em J. Biomed. Optics}, 9(3):632--647, 2004.

\bibitem{chiang}
J.~Chiang, S.~Birla, M.~Bedoya, D.~Jones, J.~Subbiah, and C.L. Brace.
\newblock Modeling and validation of microwave ablations with internal
  vaporization.
\newblock {\em IEEE Trans. Biomed. Eng.}, 62(2):657--663, 2015.

\bibitem{lc2012}
L.~Consiglieri.
\newblock Continuum models for the cooling effect of blood flow on thermal
  ablation techniques.
\newblock {\em Int. J. Thermophys.}, 33(5):864--884, 2012.

\bibitem{lc2013}
L.~Consiglieri.
\newblock An analytical solution for a bio-heat transfer problem.
\newblock {\em Int. J. Bio-Sci. Bio-Technol.}, 5(5):267--278, 2013.

\bibitem{lc2015}
L.~Consiglieri.
\newblock {Analytical solutions in the modeling of the local RF ablation}.
\newblock {\em J. Mech. Med. Biol.}, 16(05):1650071, 2016.

\bibitem{lc2024}
L.~Consiglieri.
\newblock Analytical solutions in the modelling of the endovenous laser
  ablation.
\newblock {\em Rev. Acad. Colomb. Ciencias Exactas Fis. Nat.},
  48(187):254–270, 2024.

\bibitem{consiglieri2003}
L.~Consiglieri, I.~Dos Santos, and D.~Haemmerich.
\newblock Theoretical analysis of the heat convection coefficient in large
  vessels and the significance for thermal ablative therapies.
\newblock {\em Phys. Med. Biol.}, 48(24):4125--4134, 2003.

\bibitem{durk}
J.W.~Jr. Durkee, P.P. Antich, and C.E. Lee.
\newblock {Exact solutions to the multiregion time-dependent bioheat equation.
  I: Solution development}.
\newblock {\em Phys. Med. Biol.}, 35(7):847--867, 1990.

\bibitem{dutta}
J.~Dutta and B.~Kundu.
\newblock An improved analytical model for heat flow in cancerous tumours to
  avoid thermal injuries during hyperthermia.
\newblock {\em Proc. Inst. Mech. Eng. H}, 235(5):500--514, 2021.

\bibitem{Sfarra}
G.~D’Alessandro, P.~Tavakolian, and S.~Sfarra.
\newblock A review of techniques and bio-heat transfer models supporting
  infrared thermal imaging for diagnosis of malignancy.
\newblock {\em Appl. Sci.}, 2024(14), 2024.

\bibitem{fan}
Y~Fan, L~Xu, S~Liu, J~Li, J~Xia, X~Qin, Y~Li, T~Gao, and X~Tang.
\newblock The state-of-the-art and perspectives of laser ablation for tumor
  treatment.
\newblock {\em Cyborg Bionic Syst.}, 5:0062, 2024.

\bibitem{Franco}
W.~Franco, J.~Liu, G.X. Wang, J.S. Nelson, and G.~Aguilar.
\newblock Radial and temporal variations in surface heat transfer during
  cryogen spray cooling.
\newblock {\em Phys. Med. Biol.}, 50(2):387, 2005.

\bibitem{g-suarez}
A.~González-Suárez, J.J. Pérez, R.M. Irastorza, A.~D'Avila, and E.~Berjano.
\newblock Computer modeling of radiofrequency cardiac ablation: 30 years of
  bioengineering research.
\newblock {\em Comput. Methods Programs Biomed.}, 214:106546, 2022.

\bibitem{guo-kim}
Z.~Guo and K.H. Kim.
\newblock Ultrafast laser radiation transfer in heterogeneous tissues with the
  discrete ordinate method.
\newblock {\em Appl. Opt.}, 42(16):2897--2905, 2003.

\bibitem{HeBischof}
X.~He and J.C. Bischof.
\newblock Quantification of temperature and injury response in thermal therapy
  and cryosurgery.
\newblock {\em Crit. Rev. Biomed. Eng.}, 31(5 \& 6):355--421, 2003.

\bibitem{ho}
C.-S. Ho, K.-C. Ju, T.-Y. Cheng, Y.-Y. Chen, and W.-L. Lin.
\newblock {Thermal therapy for breast tumors by using a cylindrical ultrasound
  phased array with multifocus pattern scanning: A preliminary numerical
  study}.
\newblock {\em Phys. Med. Biol.}, 52:4585--4599, 2007.

\bibitem{jacq}
S.L. Jacques.
\newblock {Optical properties of biological tissues: A review}.
\newblock {\em Phys. Med. Biol.}, 58:R37--R61, 2013.

\bibitem{jaun}
M.K. Jaunich, S.~Raje, K.~Kim, K.~Mitra, and Z.~Guo.
\newblock Bio-heat transfer analysis during short pulse laser irradiation of
  tissues.
\newblock {\em Int. J. Heat Mass Transf.}, 51:5511--5521, 2008.

\bibitem{kumar}
M.~Kumar and K.N. Rai.
\newblock Three phase bio-heat transfer model in three-dimensional space for
  multiprobe cryosurgery.
\newblock {\em J. Therm. Anal. Calorim.}, 147(24):14491--14507, 2022.

\bibitem{Laza}
A.~Lazarovich, V.~Viswanath, A.S. Dahmen, and A.~Sidana.
\newblock A narrative clinical trials review in the realm of focal therapy for
  localized prostate cancer.
\newblock {\em Transl. Cancer Res.}, 13(11):6529--6539, 2024.

\bibitem{liu-boas}
H.~Liu, D.A. Boas, Y.~Zhang, A.G. Yodh, and B.~Chance.
\newblock {Determination of optical properties and blood oxygenation in tissue
  using continuous NIR light}.
\newblock {\em Phys. Med. Biol.}, 40(11):1983--1993, 1995.

\bibitem{liu}
J.~Liu, X.~Chen, and L.X. Xu.
\newblock New thermal wave aspects on burn evaluation of skin subjected to
  instantaneous heating.
\newblock {\em IEEE Trans. Biomed. Eng.}, 46(4):420--428, 1999.

\bibitem{Loiola}
B.R. Loiola, H.R.B. Orlande, and G.S. Dulikravich.
\newblock Thermal damage during ablation of biological tissues.
\newblock {\em Numer. Heat Transf. A}, 73(10):685--701, 2018.

\bibitem{manenti}
G.~Manenti, T.~Perretta, M.~Nezzo, F.R. Fraioli, B.~Carreri, P.E. Gigliotti,
  A.~Micillo, A.~Malizia, D.~Di~Giovanni, C.P. Ryan, and F.G. Garaci.
\newblock {Transperineal Laser Ablation (TPLA) treatment of focal
  low-intermediate risk prostate cancer}.
\newblock {\em Cancers}, 16(1404), 2024.

\bibitem{marq}
M.-F. Marqa, P.~Colli, P.~Nevoux, S.R. Mordon, and N.~Betrouni.
\newblock {Focal Laser Ablation of prostate cancer: Numerical simulation of
  temperature and damage distribution}.
\newblock {\em BioMed. Eng. OnLine}, 10(45), 2011.

\bibitem{mordon}
S.~Mordon, B.~Buys, J.M. Brunetaud, and Y.~Moschetto.
\newblock New directions in medical laser concept: role of laser-tissue
  interaction modelization and feedback control.
\newblock {\em Lasers Med. Sci.}, 4(Suppl 1):317--327, 1989.

\bibitem{niemz}
M.H. Niemz.
\newblock {\em Laser-tissue interactions: Fundamentals and applications}.
\newblock Springer-Verlag, 2007.

\bibitem{oden}
J.T. Oden, K.R. Diller, C.Bajaj, J.C. Browne, J.~Hazle, I.Babu\v ska, J.~Bass,
  L.~Bidual, L.~Domkowicz, A.~Elliott, Y.~Fang, D.~Fuentes, S.~Prudhomme, M.N.
  Rylander, R.J. Stafford, and Y.~Zhang.
\newblock Dynamic data-driven finite element models for laser treatment of
  cancer.
\newblock {\em Numer. Methods Partial Differ. Equ.}, 23(4):904--922, 2007.

\bibitem{olver}
F.W.J. Olver.
\newblock {\em Bessel functions of integer order}, chapter~9, pages 355--389.
\newblock United States Department of Commerce, Washington, D.C., 1972.

\bibitem{ozisik}
M.N. \"Ozisik.
\newblock {\em Heat conduction}.
\newblock John Wiley \& Sons, Ltd, 2012.

\bibitem{paul}
A.~Paul, N.K. Bandaru, A.~Narasimhan, and S.K. Das.
\newblock Subsurface tumor ablation with near-infrared radiation using
  intratumoral and intravenous injection of nanoparticles.
\newblock {\em Int. J. Micro-Nano Scale Transport}, 5(2):69--80, 2014.

\bibitem{w-vang7}
S.A.. Prahl.
\newblock {\em The diffusion approximation in three dimensions}, chapter~7,
  pages 207--231.
\newblock Springer US, Boston, MA, 1995.

\bibitem{pupo}
A.E.B. Pupo, M.M. Gonz\'alez, L.E.B. Cabrales, J.B. Reys, E.J.R. Oria, J.J.G.
  Nava, R.P. Jim\'enez, F.M. Sanch\'ez, H.M.C. Ciria, and J.M.B. Cabrales.
\newblock 3d current density in tumors and surrounding healthy tissues
  generated by a system of straight electrode arrays.
\newblock {\em Math. Comput. Simulation}, 138:49--64, 2017.

\bibitem{putzer}
M.~Putzer, G.R. da~Silva, K.~Michael, N.~Schröder, T.~Schudeleit, M.~Bambach,
  and K.~Wegener.
\newblock Geometrical modeling of ultrashort pulse laser ablation with
  redeposition and analysis of the influence of spot size and shape.
\newblock {\em Mater. Des.}, 246:113357, 2024.

\bibitem{ryl2005}
M.N. Rylander, Y.~Feng, J.~Bass, and K.R. Diller.
\newblock Thermally induced injury and heat shock protein expression in cells
  and tissues.
\newblock {\em Ann. N.Y. Acad. Sci.}, 1066:222--242, 2005.

\bibitem{ryl}
M.N. Rylander, Y.~Feng, K.~Zimmermann, and K.R. Diller.
\newblock Measurement and mathematical modeling of thermally induced injury and
  heat shock protein expression kinetics in normal and cancerous prostate
  cells.
\newblock {\em Int. J. Hyperth.}, 26(8):748--764, 2010.

\bibitem{szhu}
J.L. Sandell and T.C. Zhu.
\newblock {A review of \textit{in-vivo} optical properties of human tissues and
  its impact on PDT}.
\newblock {\em J. Biophotonics}, 4(11-12):773--787, 2011.

\bibitem{schena}
E.~Schena, P.~Saccomandi, and Y.~Fong.
\newblock {Laser Ablation for cancer: Past, present and future}.
\newblock {\em J. Funct. Biomater.}, 8(2), 2017.

\bibitem{w-vang6}
W.M. Star.
\newblock {\em Diffusion theory of light transport}, chapter~6, pages 131--206.
\newblock Springer US, New York, 1995.

\bibitem{tang}
D.W. Tang and N.~Araki.
\newblock {On non-Fourier temperature wave and thermal relaxation time}.
\newblock {\em Int. J. Thermophys.}, 18(2):493--504, 1997.

\bibitem{DeVita}
E.~De Vita, M.~De Landro, C.~Massaroni, A.~Iadicicco, P.~Saccomandi, E.~Schena,
  and S.~Campopiano.
\newblock Fiber optic sensors-based thermal analysis of perfusion-mediated
  tissue cooling in liver undergoing laser ablation.
\newblock {\em IEEE Trans. Biomed. Eng.}, 68(3), 2021.

\bibitem{voel}
R.K. Voeller, R.B. Schuessler, and R.J.~Damiano Jr.
\newblock {\em {Surgical treatment of atrial fibrillation}}, chapter~59, pages
  1375--1394.
\newblock McGraw-Hill, New York, 2008.

\bibitem{vogel}
A.~Vogel and V.~Venugopalan.
\newblock {\em Pulsed Laser Ablation of soft biological tissues}, chapter~14,
  pages 551--615.
\newblock Springer Science+Business Media B.V., 2011.

\bibitem{WangWuWu}
J.~Wang, S.~Wu, Z.~Wu, H.~Gao, and S.~Huang.
\newblock Influences of blood flow parameters on temperature distribution
  during liver tumor microwave ablation.
\newblock {\em Front. Biosci.-Landmark}, 26(9):504--516, 2021.

\bibitem{whit1d}
P.~Whiting, J.M. Dowden, P.D. Kapadia, and M.P. Davis.
\newblock A one-dimensional mathematical model of laser induced thermal
  ablation of biological tissue.
\newblock {\em Lasers Med. Sci.}, 7:357--368, 1992.

\end{thebibliography}

\end{document}